\documentclass[final]{siamltex}
\usepackage{amssymb}
\usepackage{amsfonts}
\newtheorem{F}{\bf Condition}[section]

\title{An example on deformation of a pair}

\author{ B. Wang\\
Nov 1, 2014}

\begin{document}

\maketitle

\begin{abstract}
In this paper we give an example to show Clemens' conjecture is not a first order
deformation problem.

\end{abstract}

\pagestyle{myheadings}
\thispagestyle{plain}
\markboth{ B. Wang}{An example on deformation of a pair}

To give a rigorous explaination of the abstract, we start with Clemens' conjecture. 
   Let $X_0$ be a smooth quintic of $\mathbf P^4$ over $\mathbb C$, and $C_0$ be an irreducible rational curve of degree $d$. 
  In [1], its original 1986 statement, Clemens 
proposed:\par
(1) the generic quintic threefold $X_0$ admits only finitely many rational curves $C_0$ of each degree.\par
(2) Each rational curve $C_0$ is a smoothly embedded $\mathbf P^1$ with normal bundle
\begin{equation}
\mathcal O_{\mathbf P^1}(-1)\oplus  O_{\mathbf P^1}(-1).\end{equation}\par
(3) All the rational curves on $X_0$ are mutually disjoint. The number of rational curves of degree $d$ on $X_0$ is
\begin{equation}
(interesting \ number)\cdot 5^3\cdot d.
\end{equation}
\bigskip

In [2], we proved
\bigskip

\begin{theorem}
Let $X_0\subset \mathbf P^4$ be a generic quintic threefold. For each degree $d\geq 1$, \par
(i) there are only finitely many irreducible rational curves $C_0\subset X_0$ of degree $d$. 
\par
(ii) Each rational curve in (i) is an immersed rational curve with normal bundle
$$N_{C_0/X_0}\simeq \mathcal O_{\mathbf P^1}(-1)\oplus \mathcal O_{\mathbf P^1}(-1).
$$
By  ``immersed rational curve" we mean that the normalization map is an immersion. 
\end{theorem}
\bigskip

Our proof uses higher order deformations of the pair, i.e. we used the genericity of the quintic $X_0$. 
This paper concerns this assumption of the conjecture, i.e.  $X_0$ is generic. The question is \par
{\em In Clemens' conjecture, can the ``genericity" of the quintic $X_0$ be replaced by some other more specific condition on the pair?}
\medskip

We are far from being ready to even formulate any conjecture on the question. This paper gives a negative answer to a specific conjectural 
condition. 
Let's start with this  first order condition.
Let $H^0(\mathcal O_{\mathbf P^4}(5))$ denote the vector space of homogeneous polynomials of degree $5=deg(X_0)$ in $5$ variables. 
Let  $f_0\in H^0(\mathcal O_{\mathbf P^4}(5))$ such that $$X_0=div(f_0)$$
is a smooth hypersurface. Let $$[f_0] \in S=\mathbf P(H^0(\mathcal O_{\mathbf P^4}(5)))$$
denote the corresponding point of $f_0$ in the projectivization.  Let 
\begin{equation}
M_d= H^0(\mathcal O_{\mathbf P^1}(d))^{\oplus 5}.\end{equation}
So the normalization of $C_0$, $c_0: \mathbf P^1\to C_0$ lies in the projectivization of $M_d$. 
Let  \begin{equation}
\Gamma \end{equation}
be the incidence scheme
\begin{equation}
\{ (c, [f]): f(c(t))=0\}\subset M_d\times S
\end{equation}
containing the point $(c_0, [f_0])$. 
\bigskip

\begin{F}: 
The projection 
\begin{equation}\begin{array}{ccc}
T_{(c_0, f_0)}\Gamma &\stackrel{\pi} \rightarrow & T_{[f_0]}S
\end{array}\end{equation} is 
surjective. 
\end{F}

\bigskip

Motivation for the condition 0.1: \par
(1) The main motivation is the question above. Can the genericity of the quintic be replaced by condition 0.1 ? There is an evidence for an affirmative answer.  Condition 0.1 was used in [3], [4]. In particular it is shown in [3] that the condition 0.1 is responsible for 
the global generation of the twisted normal bundle $N_{C_0/X_0}(d)$. 
So it is natural to conjecture it may also be responsible for the 
vanishing of  $H^1$ of the un-twisted normal bundle $N_{C_0/X_0}$.
\par

(3) If $X_0$ is generic, the condition 0.1 is automatic. Thus the studying only the first orders 
in Clemens' conjecture is equivalent to assuming the condition 0.1 and removing the genericity condition of
the quintic.

\bigskip

However the following theorem shows that expectations above are wrong.\bigskip

 \begin{theorem} 
There are a smooth quintic threefold $X_0$ and a plane rational curve $C_0\subset X_0$ of degree $3$ with
the normalization $c_0: \mathbf P^1\to C_0$, such that
the  condition 0.1 is satisfied, i.e. 
the projection 
\begin{equation}\begin{array}{ccc}
T_{(c_0, f_0)}\Gamma &\stackrel{\pi} \rightarrow & T_{[f_0]}S
\end{array}\end{equation} is 
surjective, but $H^1(N_{c_0/X_0})\neq 0$, where
$N_{c_0/X_0}$ is defined to be the normal sheaf of the normalization map $c_0$. 
\end{theorem}

\bigskip

{\bf Remark}  Actually it is not so difficult to see that for a  hypersurface $X_0$ of another types (such that Fano) the condition 0.1 is insufficient for the
vanishing of $H^1$ of the normal sheaf. But  this becomes non-trivial for a quintic. It was proved in [3] that the condition 0.1  indeed implies the Clemens' conjecture for lines. We believe that for rational curves of degree less than 3, the condition 0.1 should be sufficient for Clemens' conjecture. 
Thus $3$ is expected to be the minimum degree for the condition 0.1 to be insufficient for Clemens' conjecture.
\bigskip

\begin{proof}
Since $c_0: \mathbf P^1\to \mathbf P^4$ generically is an embedding, the normal sheaves 
\begin{equation}
N_{c_0/\mathbf P^4}, N_{c_0/X_0}
\end{equation}
are well-defined and there is an exact seqence of sheaves over $\mathbf P^1$, 
\begin{equation}\begin{array}{ccccccccc}
0 &\rightarrow & N_{c_0/X_0} &\rightarrow &N_{c_0/\mathbf P^4} &\rightarrow &
c_0^\ast( N_{X_0/\mathbf P^4})&\rightarrow & 0.
\end{array}\end{equation}
This induces an exact sequence of vector spaces, 

\begin{equation}\begin{array}{ccccccc}
 H^0(N_{c_0/\mathbf P^4}) &\stackrel{\nu^s} \rightarrow &H^0(c_0^\ast( N_{X_0/\mathbf P^4})) &\rightarrow &
H^1 (N_{c_0/X_0}) &\rightarrow & 0.
\end{array}\end{equation}
We need to add a few more words about $\nu^s$ for the calculation later.
In [2], we proved that if $c_0$ is a birational to its image, there is a natural surjection
\begin{equation}\begin{array}{ccc}
H^0(\mathcal O_{\mathbf P^1}(d))^{\oplus 5} &\rightarrow & H^0(c_0^\ast (T_{\mathbf P^1}))
\end{array}
\end{equation}
whose kernel is a one-dimensional line through the origin.  Hence we obtain a composition 
\begin{equation}\begin{array}{ccc}
\phi: H^0(\mathcal O_{\mathbf P^1}(d))^{\oplus 5} &\stackrel{\phi}\rightarrow &
H^0(c_0^\ast( N_{X_0/\mathbf P^4}))\\
(\alpha_0, \alpha_4) &\rightarrow & \sum_{i=0}^4 {\partial c_0^\ast (f_0)\over \partial \alpha_i}.
\end{array}\end{equation}
(The partial derivative ${\partial c_0^\ast ( f_0)\over \partial \alpha_i}$ denotes $ {\partial c^\ast ( f_0)\over \partial \alpha_i}|_{c_0}$).
In (0.12), we identify \begin{equation}
T_0 (H^0 (\mathcal O_{\mathbf P^1}(d)))
\end{equation}
with 
\begin{equation}
H^0 (\mathcal O_{\mathbf P^1}(d)). 
\end{equation}
by the linearity of $H^0 (\mathcal O_{\mathbf P^1}(d))$. 
Furthermore $\nu^s$ is surjective if and only if $\phi$ is surjective which is equivalent to
the kernel of $\phi$ being the tangent space $T_I(GL(2))$ of the orbit by the  $GL(2)$ action on
$$H^0(\mathcal O_{\mathbf P^1}(d))^{\oplus 5}$$
induced from the automorphsims of $\mathbb C^2$ for $\mathbf P(\mathbb C^2)=\mathbf P^1$. 

\bigskip

On the other hand, there is a linear map $\mu$ at $c_0$, defined by
\begin{equation}\begin{array}{ccc}
H^0(\mathcal O_{\mathbf P^4}(5)) &\stackrel{\mu}\rightarrow & H^0( c_0^\ast (N_{X_0/\mathbf P^4}))\\
f &\rightarrow &c_0^\ast(f).\end{array}\end{equation}
The differential $\mu^s$ of it is also a linear map
 $$\begin{array}{ccc}
T_{f_0}H^0(\mathcal O_{\mathbf P^4}(5))=H^0(\mathcal O_{\mathbf P^4}(5)) &\stackrel{\mu^s}\rightarrow & T_0 H^0( c_0^\ast (N_{X_0/\mathbf P^4}))
=H^0( c_0^\ast (N_{X_0/\mathbf P^4}))\\
\alpha &\rightarrow  & {\partial c_0^\ast(f_0)\over \partial \alpha} .\end{array}$$

Because of (0.10), $H^1 (N_{c_0/X_0})=0$ is equivalent to the surjectivity of $\nu^s$. At the meantime, the condition 0.1
can be expressed as
\begin{equation}
image(\mu^s)\subset image(\nu^s).
\end{equation}

Therefore we would like to construct an example such that
(0.16) holds but $\nu^s$ is not surjective.
\bigskip

The following is this example. 
 Let 
$z_0, \cdots, z_4$ be the homogeneous coordinates of $\mathbf P^4$. 
Let $V\subset \mathbf P^4$ be the 
plane \begin{equation}
\{[z_0\cdots, z_4]: z_0=z_1=0\}.
\end{equation}
There is a natural embedding
\begin{equation}
H^0(\mathcal O_V(n))\stackrel{\j}\hookrightarrow H^0(\mathcal O_{\mathbf P^4}(n)).
\end{equation}

Let $C_0$ be a generic cubic rational curve on $V$, i.e. $C_0$ is defined by a polynomial
\begin{equation}
g_2\in H^0(\mathcal O_V(3))
\end{equation}
that has exactly one nodal singularity. 
Let $c_0:\mathbf P^1\rightarrow C_0$ be its normalization.
Let $g_0\in H^0(\mathcal O_{\mathbf P^4}(5))$, $q\in H^0(\mathcal O_V(2))$ such that
\begin{equation}
\xi ^\ast (g_0)
\end{equation}
 and $c_0^\ast( q)$ have at least three common zeros, where
$\xi$  is the embedding $$\mathbf P^1\hookrightarrow C_0\hookrightarrow V\hookrightarrow \mathbf P^4.$$

Let
\begin{equation}\begin{array}{c} 
f_0= z_0 g_0 +z_1 g_1+g_2 \j(q)

\end{array}\end{equation}
where $g_1\in H^0(\mathcal O_{\mathbf P^4}(5))$ is generic, and all other polynomials 
$g_0, g_2, q$ are generic in the chosen types above.

Then we can check that \par
(1) $X_0=div(f_0)$ is a smooth quintic threefold, \par
(2) $C_0=im(c_0)$ is a singular, cubic plane rational curve.  But it is  a global immersion, \par
(3) $C_0$ lies on $X_0$.

\smallskip

In the following calculation, we identify 
\begin{equation}
T_0 (H^0 (\mathcal O_{\mathbf P^1}(m)))
\end{equation}
with 
\begin{equation}
H^0 (\mathcal O_{\mathbf P^1}(m)). 
\end{equation}
by the linearity of $H^0 (\mathcal O_{\mathbf P^1}(m))$. 

Consider the linear map $\lambda$

\begin{equation}\begin{array}{ccc}
 H^0(c_0^\ast(\mathcal O_{\mathbf V}(1)))^{\oplus 5} &\stackrel{\lambda} \rightarrow &  H^0( c_0^\ast (\mathcal O_{V}(5)))\\
c_0^\ast (\alpha_0, \cdots, \alpha_4)  &\rightarrow & c_0^\ast ( \sum_{i=0}^1 {\partial g_i|_{V}\over \partial \alpha_i}+q\sum_{i=2}^4 {\partial g_2\over \partial \alpha_i}  ).
\end{array}\end{equation}

Because of the genericity of $g_0, g_1, g_2$, the linear map $\lambda$ is surjective. This shows that
the condition 0.1 is satisfied. 

Because \begin{equation}
\xi ^\ast (f_0)
\end{equation}
 and $c_0^\ast(q)$ have at least three common zeros, we find non-zero $$\beta_0, \beta_2, \beta_3, \beta_4 \in H^0 (\mathcal O_{\mathbf P^1}(3))$$
such that $$(\beta_0, 0, \beta_2, \beta_3, \beta_4)$$ is not in the tangent space $T_I(GL(2))$ of the automorphism $GL(2)$ of
\begin{equation}
H^0 (\mathcal O_{\mathbf P^1}(3))^{\oplus 5}
\end{equation}
(because $\beta_0$ is not identically zero)
and
\begin{equation}
\beta_0  \cdot c_0^\ast( g_0)+ \xi ^\ast (q) \cdot \sum_{i=2}^4 \beta_i\cdot {\partial g_2\over \partial z_i}  =0.\end{equation}
Hence $$(\beta_0, 0, \beta_2, \beta_3, \beta_4)$$ is in the kernel of $\phi$ but not in the $T_I(GL(2))$. 
Then \begin{equation} dim(ker(\phi))\geq 5.
\end{equation} 
This shows 
$\nu^s$ for this pair $X_0, C_0$ is not surjective.

The theorem is proved.

\end{proof} 

\bigskip

This example  offers a scenario that \begin{equation}
image(\mu^s)=image(\nu^s)\subsetneq H^0(\mathcal O_{\mathbf P^1}(15)).
\end{equation}

\end{document}